\documentclass[12pt]{amsart}

\usepackage{amssymb,euscript,amsmath, mathrsfs,amscd}
\usepackage[dvips]{graphicx}
\usepackage{epsfig,pstricks}
\usepackage{fullpage}

\newcounter{ENUM}
\newcommand{\be}{\begin{enumerate}}
\newcommand{\ee}{\end{enumerate}}
\newcommand{\beas}{\begin{eqnarray*}}
\newcommand{\eeas}{\end{eqnarray*}}
\newcommand{\beq}{\begin{equation}}
\newcommand{\eeq}{\end{equation}}

\newcommand{\st}{\,:\,}

\newtheorem{thm}{Theorem}[section]

\newtheorem{lem}[thm]{Lemma}

\theoremstyle{definition}

\newtheorem{ex}[thm]{Example}

\theoremstyle{remark}

\numberwithin{equation}{section}

\def\zz{\mathbb{Z}}

\newcommand{\jt}{\mathrm{JT}}
\newcommand{\jtl}{\mathrm{JT}_\lambda}
\newcommand{\qq}{\mathbb{Q}}
\newcommand{\vjt}{\varphi_n\mathrm{JT}_\lambda}
\newcommand{\jtq}{\mathrm{JT}(q)_\lambda}


\begin{document}
\title[Smith Normal Form]{The Smith Normal Form of a Specialized
Jacobi-Trudi Matrix}

\date{\today}

\author{Richard P. Stanley}
\email{rstan@math.mit.edu}
\address{Department of Mathematics, University of Miami, Coral Gables,
FL 33124}

\thanks{Partially supported by NSF grant DMS-1068625.}

\begin{abstract}
Let JT$_\lambda$ be the Jacobi-Trudi matrix corresponding to the
partition $\lambda$, so $\det\jt_\lambda$ is the Schur function
$s_\lambda$ in the variables $x_1,x_2,\dots$. Set $x_1=\cdots=x_n=1$
and all other $x_i=0$. Then the entries of $\jt_\lambda$ become
polynomials in $n$ of the form $\binom{n+j-1}{j}$. We determine the
Smith normal form over the ring $\qq[n]$ of this specialization of
$\jtl$ . The proof carries over to the specialization $x_i=q^{i-1}$
for $1\leq i\leq n$ and $x_i=0$ for $i>n$, where we set $q^n=y$ and
work over the ring $\qq(q)[y]$.   
\end{abstract}

\maketitle

\section{Introduction}

Let $M$ be an $r\times s$ matrix over a commutative ring $R$ (with
identity), which for convenience we assume has full rank $r$. If there
exist invertible $r \times r$ and $s \times s$ matrices $P$ and $Q$
such that the product $PMQ$ is a diagonal matrix with diagonal entries
$\alpha_1,\alpha_2,\dots, \alpha_r$ satisfying $\alpha_i \mid
\alpha_{i+1} $ for all $1 \le i \le r-1$, then $PMQ$ is called the
\emph{Smith normal form (SNF)} of $M$.  In general, the SNF does not
exist. It does exist when $R$ is a principal ideal domain (PID) such
as $\qq[n]$, the polynomial ring in the indeterminate $n$ over the
rationals (which is the case considered in this paper). Over a PID the
SNF is unique up to multiplication of diagonal elements by units in
$R$. Note that the units of the ring $\qq[n]$ are the nonzero rational
numbers. Since the determinants of $P$ and $Q$ are units in $R$, we
obtain when $M$ is a nonsingular square matrix a canonical
factorization $\det M = u\alpha_1\alpha_2\cdots \alpha_m$, where $u$
is a unit. Thus whenever $\det M$ has a lot of factors, it suggests
that it might be interesting to consider the SNF.

There has been a lot of recent work, such as \cite{b-s}\cite{w-s}, on
the Smith normal form of specific matrices and random matrices, and on
different situations in which SNF occurs. Here we will determine the
SNF of a certain matrices that arise naturally in the theory of
symmetric functions. We will follow notation and terminology from
\cite[Chap.~7]{ec2}. Namely, let $\lambda=(\lambda_1,\lambda_2,\dots)$
be a partition of some positive integer, and let $h_i$ denote the
complete homogeneous symmetric function of degree $i$ in the variables
$x_1, x_2,\dots$. Set $h_0=1$ and $h_m=0$ for $m<0$. Let $t$ be an
integer for which $\ell(\lambda)\leq t$, where $\ell(\lambda)$ denotes
the length (number of parts) of $\lambda$. The \emph{Jacobi-Trudi
  matrix} $\jtl$ is defined by
  $$ \jtl=\left[ h_{\lambda_i+j-i}\right]_{i,j=1}^t. $$
The \emph{Jacobi-Trudi identity} \cite[{\textsection}7.16]{ec2}
asserts that 
$\det \jtl =s_\lambda$, the Schur function indexed by $\lambda$.

For a symmetric function $f$, let $\varphi_n f$ denote the
specialization $f(1^n)$, that is, set $x_1=\cdots=x_n=1$ and all other
$x_i=0$ in $f$. It is easy to see \cite[Prop.~7.8.3]{ec2} that
   \beq \varphi_n h_i =\binom{n+i-1}{i}, \label{eq:phihi} \eeq
a polynomial in $n$ of degree $i$. Identify $\lambda$ with its (Young)
diagram, so the squares of $\lambda$ are indexed by pairs $(i,j)$,
$1\leq i\leq \ell(\lambda)$, $1\leq j\leq \lambda_i$. The
\emph{content} $c(u)$ of the square $u=(i,j)$ is defined to be
$c(u)=j-i$. A standard result \cite[Cor.~7.21.4]{ec2} in the theory of
symmetric functions states that
  \beq \varphi_n s_\lambda =\frac{1}{H_\lambda}\prod_{u\in\lambda}
   (n+c(u)), \label{eq:hc} \eeq
where $H_\lambda$ is a positive integer whose value is irrelevant here
(since it is a unit in $\qq[n]$). Since this polynomial factors a lot
(in fact, into linear factors) over $\qq[n]$, we are motivated to
consider the SNF of the matrix 
  $$ \varphi_n\jtl = \left[ \binom{n+\lambda_i+j-i-1}{\lambda_i+j-i}
     \right]_{i,j=1}^t. $$

Let $D_k$ denote the $k$th \emph{diagonal hook} of $\lambda$, i.e.,
all squares $(i,j)\in\lambda$ such that either $i=k$ and $j\geq k$, or
$j=k$ and $i\geq k$. Note that $\lambda$ is a disjoint union of its
diagonal hooks. If $r=\mathrm{rank}(\lambda)\mathrel{\mathop:}= \max\{
i\st \lambda_i\geq i\}$, then note also that $D_k=\emptyset$ for
$k>r$. Our main result is the following.

\begin{thm} \label{thm1} 
Let the SNF of $\varphi_n\jtl$ have main diagonal
$(\alpha_1,\alpha_2,\dots,\alpha_t)$, where $t\geq
\ell(\lambda)$. Then we can take
  $$ \alpha_i = \prod_{u\in D_{t-i+1}} (n+c(u)). $$
\end{thm}

An equivalent statement to Theorem~\ref{thm1} is that the $\alpha_i$'s
are squarefree (as polynomials in $n$), since $\alpha_t$ is the
largest squarefree factor of $\varphi_n s_\lambda$, $\alpha_{t-1}$ is
the largest squarefree factor of $(\varphi_n s_\lambda)/\alpha_t$, etc.

\begin{ex} 
Let $\lambda=(7,5,5,2)$. Figure~\ref{fig1} shows the diagram of
$\lambda$ with the content of each square. Let $t=\ell(\lambda)=4$. We
see that 
  \beas \alpha_4 & = & (n-3)(n-2)\cdots (n+6)\\
        \alpha_3 & = & (n-2)(n-1)n(n+1)(n+2)(n+3)\\
        \alpha_2 & = & n(n+1)(n+2)\\
        \alpha_1 & = & 1. \eeas

\begin{figure}
\centering
\centerline{\includegraphics[width=6cm]{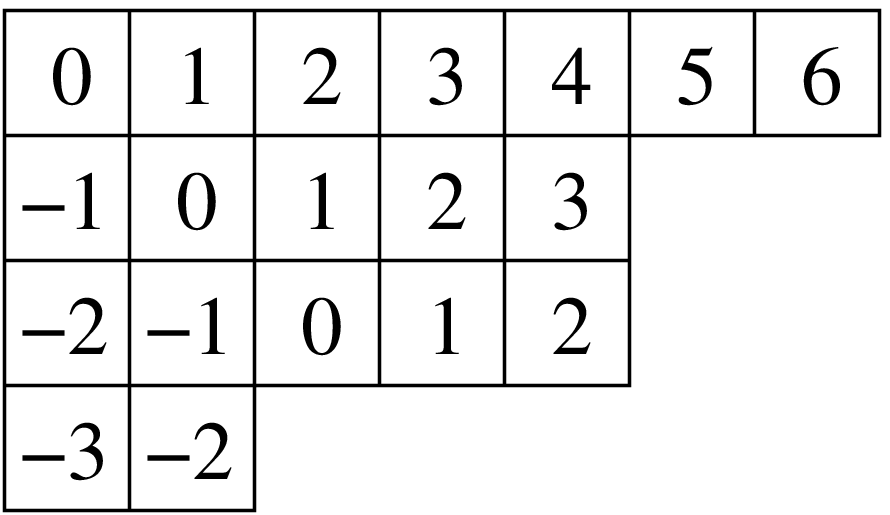}}
\caption{The contents of the partition $(7,5,5,2)$}
\label{fig1}
\end{figure}

\end{ex}

The problem of computing the SNF of a suitably specialized
Jacobi-Trudi matrix was raised by Kuperberg \cite{kup}. His Theorem~14
has some overlap with our Theorem~\ref{thm1}. Propp
\cite[Problem~5]{propp} mentions a two-part question of Kuperberg. The
first part is equivalent to our Theorem~\ref{thm1} for rectangular
shapes. (The second part asks for an interpretation in terms of
tilings, which we do not consider.)

\section{Proof of the main theorem}
To prove Theorem~\ref{thm1} we use the following well-known
description of SNF over a PID.

\begin{lem} \label{lem:gcd}
Let $\mathrm{diag}(\alpha_1,\dots, \alpha_m)$ be the SNF of an
$m\times n$ matrix $M$ over a PID.
Then $\alpha_1 \alpha_2 \cdots \alpha_k$ is the
greatest common divisor (gcd) of the $k\times k$ minors of $M$. 
\end{lem}

Let $\lambda$ be a partition of length at most $t$ and with diagonal
hooks $D_1,\dots, D_t$.  Given the $t\times t$ matrix $\varphi_n\jtl$
and $1\leq k\leq t$, let $M_k$ be the square submatrix consisting of
the last $k$ rows and first $k$ columns of $\varphi_n\jtl$.  We claim
the following.

\be
 \item[C1.] If $\det M_k=0$ then $\vjt$ has a $k\times k$ minor equal
   to 1.\ \ Otherwise,
   \beq \det M_k = c_k\prod_{i=1}^k \prod_{i\in D_{t-i+1}} (n+c(u)), 
     \label{eq:detmk} \eeq
where $c_k$ is a nonzero rational number.  
 \item[C2.] If $\det M_k\neq 0$, then every $k\times k$ minor of
   $\vjt$ is divisible (in the ring $\qq[n]$) by $\det M_k$.
\ee

\emph{Proof of C1.} It is well known and follows immediately from the
Jacobi-Trudi identity for skew Schur functions that every minor of
$\jtl$ is either 0 or a skew Schur function $s_{\rho/\sigma}$ for some
skew shape $\rho/\sigma$. Let $N$ be a $k\times k$ submatrix of $\jtl$
with determinant zero. This can only happen if $N$ is strictly upper
triangular, since otherwise the determinant is a nonzero
$s_{\rho/\sigma}$. Each row of $\jtl$ that intersects $N$ consists of
a string of 0's, followed by a 1, and possibly followed by other
terms. The 1's in these rows appear strictly from left-to-right as we
move down $\jtl$. Hence the $k\times k$ submatrix of $\jtl$ with the
same rows as $N$ and with each column containing 1 is upper
unitriangular and hence has determinant 1.\ \ Since $\varphi_n
s_{\rho/\sigma}\neq 0$, the same reasoning applies to $\vjt$, so the
first assertion of (C1) is proved.

If on the other hand $\det M_k\neq 0$, then $M_k$ is just the
Jacobi-Trudi matrix for the subshape $\bigcup_{i=1}^k D_{t-i+1}$ of
$\lambda$, so (C1) follows from equation~\eqref{eq:hc}.

\emph{Proof of C2.} Suppose that $\det M_k\neq 0$. Thus $M_k$ is the
Jacobi-Trudi matrix for the partition $\mu=\bigcup_{i=1}^k
D_{t-i+1}$. It is easy to check that any $k\times k$ submatrix of
$\jtl$ is the Jacobi-Trudi matrix of a skew shape $\rho/\sigma$ such
that (the diagram of) $\rho/\sigma$ has the following property:

  (P) There is a subdiagram $\nu$ (of an ordinary partition) of
$\rho/\sigma$ containing $\mu$, and all other squares of $\rho/\sigma$
are to the left of $\nu$.

Suppose now that $\langle s_{\rho/\sigma},s_\tau\rangle\neq 0$. We
claim that $\mu\subseteq \tau$. This will complete the proof, since
then $\det M_k = H_\mu^{-1} \prod_{u\in \mu} (n+c(u))$, and the
contents of $\mu$ form a submultiset of the contents of $\tau$.

The statement that $\langle s_{\rho/\sigma},s_\tau\rangle\neq 0$ is
equivalent to $c^{\rho}_{\sigma\tau}\neq 0$, where
$c^{\rho}_{\sigma\tau}$ is a Littlewood-Richardson coefficient
\cite[eqn.~(7.64)]{ec2}. By the Littlewood-Richardson rule as
formulated e.g.\ in \cite[Thm.~A1.3.3]{ec2}, $c^{\rho}_{\sigma\tau}$
is the number of semistandard Young tableaux (SSYT) of shape
$\rho/\sigma$ and content $\tau$ whose reverse reading word is a
lattice permutation. By Property~(P) such an SSYT must have the last
$\mu_i$ entries in row $i$ equal to $i$. Hence $\tau_i\geq \mu_i$ for
all $i$, as desired. This completes the proof of (C2).

\begin{figure}
\centering
\centerline{\includegraphics[width=6cm]{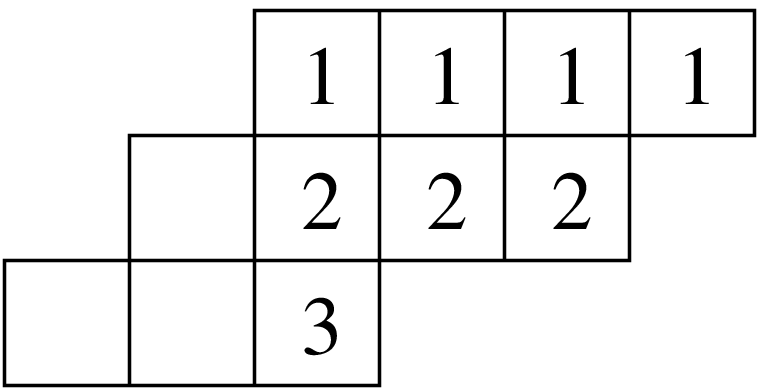}}
\caption{A partial Littlewood-Richardson filling}
\label{fig2}
\end{figure}

As an illustration of the proof of (C2), suppose that $\lambda =
(7,6,6,5,3)$ and we take $k=3$. Then $\mu=(4,3,1)$. The
$3\times 3$ minor with rows 3,4,5 and columns 1,3,5 (say) is given by
  $$ \left[ \begin{array}{ccc} h_4 & h_6 & h_8\\
      h_2 & h_4 & h_6\\ 0 & h_1 & h_3 \end{array} \right], $$
which is the Jacobi-Trudi matrix for the skew shape
$(6,5,3)/(2,1)$. Any Littlewood-Richardson filling of this shape has
to have the entries indicated in Figure~\ref{fig2}, so the type $\tau$
of this filling satisfies $\tau\supseteq (4,3,1)=\mu$.

\emph{Proof of Theorem~\ref{thm1}.} If the $k$th diagonal hook is
empty, then (C1) shows that $\jtl$ contains a $k\times k$ minor equal
to 1.\ \ Hence the gcd of the $k\times k$ minors is also 1, and
therefore the gcd of the $j\times j$ minors for each $j<k$ is 1. Thus
by Lemma~\ref{lem:gcd}, we have $\alpha_k=1$ as desired.

If the $k$th diagonal hook is nonempty, then (C2) shows that every
$k\times k$ minor is divisible by $\det M_k$. Hence the gcd of the
$k\times k$ minors is equal to $\det M_k$, and the proof follows from
equation~\eqref{eq:detmk} and Lemma~\ref{lem:gcd}. \qed

\section{A $q$-analogue}
There is a standard $q$-analogue $\varphi_n(q) s_\lambda$ of
$\varphi_n s_\lambda$ \cite[Thm.~7.21.2]{ec2}, namely, 
 \beas \varphi_n(q)s_\lambda & = & s_\lambda(1,q,q^2,\dots,q^{n-1})\\ 
  & = &
 \frac{q^{b(\lambda)}}{H_\lambda(q)}\prod_{u\in\lambda}(1-q^{n+c(u)}), \eeas 
where $H_\lambda(q)$ is a polynomial in $q$ (the $q$-analogue of
$H_\lambda$) and $b(\lambda)$ is a nonnegative integer. What is the
SNF of $\varphi_n(q)\jtl$? The problem 
arises of choosing the ring over which we compute the SNF. The most
natural choice might seem to be to fix $n$ and then work over the ring
$\qq[q]$ (or even $\zz[q]$, assuming that the SNF exists).  This
question, however, is not really a $q$-analogue of what was done
above, since we considered $n$ to be \emph{variable} while here it is
a constant. In fact, it seems quite difficult to compute the SNF this
way. Its form seems to depend on $n$ is in a very delicate
way. Instead we can set $y=q^n$. For instance,
 \beas \varphi_n(q)h_3 & = &
 \frac{(1-q^{n+2})(1-q^{n+1})(1-q^n)} {(1-q^3)(1-q^2)(1-q)}\\ & = &
 \frac{(1-q^2y)(1-qy)(1-y)}{(1-q^3)(1-q^2)(1-q)}. \eeas 
Since the entries of $\varphi_n(q)\jtl$ become polynomials in $y$ with
coefficients in the field $F=\qq(q)$, we can ask for the SNF over the
PID $F[y]$.  The proof of Theorem~\ref{thm1} carries over,
mutatis mudandi, to this $q$-version.

\begin{thm} \label{thm2}
Let $M_\lambda$ denote the matrix obtained from
$\varphi_n(q) \jtl$ by substituting $q^n=y$. Let the SNF of
$M_\lambda$ over the ring $\qq(q)[y]$ have main diagonal
$(\beta_1,\beta_2,\dots,\beta_t)$, where $t\geq\ell(\lambda)$. Then we
can take 
 $$ \beta_i = \prod_{u\in D_{t-i+1}} (1-q^{c(u)}y). $$
\end{thm}

Perhaps this result still seems to be an unsatisfactory
$q$-analogue (or in this case, a $y$-analogue) since we cannot
substitute $y=1$ to reduce to $\varphi_n
\jtl$. Instead, however, make the substitution 
  \beq y \to \frac{1}{(1-q)y+1}. \label{eq:subs} \eeq 
For any $k\in\zz$ write $\boldsymbol{(k)}=(1-q^k)/(1-q)$. For
instance, $\boldsymbol{(-3)}= -q^{-1}-q^{-2}-q^{-3}$ and
$\boldsymbol{(0)}=0$.
Under the substitution \eqref{eq:subs} we 
have for any $k\in\zz$, 
  $$ 1-q^ky\to \frac{(1-q)(y+\boldsymbol{(k)})}{(1-q)y+1}. $$
For any symmetric function $f$ let $\varphi^*f$ denote the
substitution $q^n\to 1/(1-q)y+1$ after writing $f(1,q,\dots,q^{n-1})$
as a polynomial in $q$ and $q^n$. Let $A$ be a square submatix of
$\jtl$. Since $\det A$ is a homogeneous symmetric function,  say of
degree $d$, the specialization $\varphi^*\det M$ will equal
$\left( \frac{1-q}{(1-q)y+1}\right)^d$ times the result of substituting
  \beq q^ky-1\to y+\boldsymbol{(k)} \label{eq:subs2} \eeq
in $M$ and then taking the determinant. It follows that the proof of
Theorem~\ref{thm1} also carries over for the substitution
\eqref{eq:subs2}.
We obtain the following variant of Theorem~\ref{thm2}, which is
clearly a satisfactory $q$-analogue of Theorem~\ref{thm1}.

\begin{thm} \label{thm2v}
For $k\geq 1$ let
  $$ f(k) =\frac{y(y+\boldsymbol{(1)})(y+\boldsymbol{(2)})\cdots 
     (y+\boldsymbol{(k-1)})}{\boldsymbol{(1)}\boldsymbol{(2)}
    \cdots \boldsymbol{(k)}}. $$
Set $f(0)=1$ and $f(k)=0$ for $k<0$. Define 
  $$ \jtq=\left[f(\lambda_i-i+j)\right]_{i,j=1}^t, $$
where $\ell(\lambda)\leq t$. Let the SNF of $\jtq$ over the ring
$\qq(q)[y]$ have main
diagonal $(\gamma_1,\gamma_2,\dots,\gamma_t)$. Then we can take
  $$ \gamma_i = \prod_{u\in D_{t-i+1}}(y+\boldsymbol{c(u)}). $$
\end{thm}

\end{document}